\documentclass[12pt,fleqn]{article}
\usepackage{amssymb}
\newtheorem{theorem}{Theorem}

\newcommand{\be}{\begin{equation}}
\newcommand{\ee}{\end{equation}}
\newcommand{\bea}{\begin{eqnarray}}
\newcommand{\eea}{\end{eqnarray}}
\newcommand{\bean}{\begin{eqnarray*}}
\newcommand{\eean}{\end{eqnarray*}}
\newcommand{\benu}{\begin{enumerate}}
\newcommand{\eenu}{\end{enumerate}}
\newcommand{\edo}{\end{document}}
\newcommand{\vab}{\vspace{1.3ex}}
\newcommand{\hab}{\hspace{0.2em}}

\newcommand{\mk}{\medskip}

\newcommand{\lb}{\linebreak}
\newcommand{\scst}{\scriptstyle}

\newcommand{\scsi}{\scriptsize}
\newcommand{\barr}{\begin{array}}
\newcommand{\earr}{\end{array}}
\newcommand{\gG}{\Gamma}
\newcommand{\gD}{\Delta}
\newcommand{\gF}{\Phi}
\newcommand{\gL}{\Lambda}
\newcommand{\ga}{\alpha}
\newcommand{\gb}{\beta}
\renewcommand{\gg}{\gamma}
\newcommand{\gd}{\delta}
\newcommand{\gep}{\varepsilon}
\newcommand{\gk}{\kappa}
\newcommand{\gl}{\lambda}
\newcommand{\gm}{\mu}
\newcommand{\gn}{\nu}
\newcommand{\gf}{\varphi}
\newcommand{\gr}{\rho}
\newcommand{\gs}{\sigma}
\newcommand{\gt}{\tau}
\newcommand{\go}{\omega}
\newcommand{\cC}{{\cal C}}
\newcommand{\fP}{{\mathfrak P}}
\newcommand{\fg}{\mathfrak g}
\newcommand{\fh}{\mathfrak h}
\newcommand{\fa}{\mathfrak a}
\newcommand{\fsl}{{\mathfrak{sl}}}
\newcommand{\fso}{{\mathfrak{so}}}
\newcommand{\fsp}{{\mathfrak{sp}}}
\newcommand{\fspin}{{\mathfrak{spin}}}
\newcommand{\fsu}{{\mathfrak{su}}}
\newcommand{\fu}{{\mathfrak{u}}}
\newcommand{\fpin}{{\mathfrak{pin}}}
\newcommand{\CC}{{\mathbb C}}
\newcommand{\HH}{{\mathbb H}}
\newcommand{\KK}{{\mathbb K}}
\newcommand{\RR}{{\mathbb R}}
\newcommand{\Hu}{\tilde{H}_{\tilde{u}}}
\newcommand{\VH}{V_{\tilde{H}}}
\newcommand{\tA}{\tilde{A}}

\newcommand{\tH}{\tilde{H}}

\newcommand{\tT}{\tilde{T}}
\newcommand{\tX}{\tilde{X}}

\newcommand{\tfh}{\tilde{\fh}}
\newcommand{\tu}{\tilde{u}}

\begin{document}
\newcommand{\Ad}{{\rm Ad}}
\renewcommand{\O}{{\rm O}}
\renewcommand{\Re}{{\rm Re}}
\newcommand{\SL}{{\rm SL}}
\newcommand{\SU}{{\rm SU}}
\newcommand{\SO}{{\rm SO}}
\newcommand{\so}{{\rm so}}
\newcommand{\Sp}{{\rm Sp}}
\newcommand{\Spin}{{\rm Spin}}
\newcommand{\G}{{\rm G}}
\newcommand{\GL}{{\rm GL}}
\newcommand{\U}{{\rm U}}
\newcommand{\diag}{{\rm diag}}

\title{Parallel spinors and holonomy groups on pseudo-Riemannian spin 
manifolds}
 
\author{Helga Baum and Ines Kath}
\date{}
\maketitle

\begin{abstract}
We describe the possible holonomy groups of simply connected irreducible 
non-locally symmetric pseudo-Riemannian spin manifolds which admit parallel spinors. 

\end{abstract}

\section{Introduction}\label{s-introduction}
In \cite{Wa 89} McKenzie~Y.~Wang described the possible holonomy groups of 
complete simply connected irreducible non-flat Riemannian spin manifolds $M^n$ which 
admit parallel spinors. Parallel spinors occur exactly for the holonomy 
representations $\SU(m)$, $n=2m\ge 4$, $\; \Sp(m)$, $n=4m\ge8$, $\G_2$ and $\Spin(7)$. A 
complete, simply connected locally symmetric Riemannian spin manifold with non-trivial 
parallel spinors is flat. 

In the present paper we consider the same problem for pseudo-Riemannian spin 
manifolds. We describe the possible holonomy groups of simply connected irreducible 
non-locally symmetric pseudo-Riemannian spin manifolds $M^{r,s}$ of dimension $r+s$ and 
index $r$ and study the chiral type and the causal type of these parallel spinors with 
respect to the indefinite scalar product on the spinor bundle. 

The method to attack this problem is based on the fact that the space of parallel 
spinor fields can be identified with the vector space 
\[
V_{\Hu}=\{v\in\gD_{r,s}\;|\; \gr(\Hu)v=v\}
\]
of all elements of the spinor modulee $\gD_{r,s}$ which are invariant under the 
action of the holonomy group $\Hu\subset \Spin(r,s)$ of the Levi-Civita connection 
$\tilde{A}$ on the spin structure $Q$ of $M^{r,s}$ with respect to $\tilde{u}\in Q$. 
The
first
list of all possible holonomy groups $H_u$ of simply connected irreducible 
non-locally symmetric semi-Riemannian manifolds was established by M.~Berger in \cite{Be 
55}. 
Later on there were made some corrections and additions to Berger's original list 
by 
D.~Alekseevskii~\cite{A 68}, R.~Brown and A.~Gray~\cite{BG 72}, R.~McLean and 
by R.~Bryant~\cite{Br 96}. 
Since $\gl(\Hu)=H_{\gl(\tilde{u})}$, where $\gl\,:\, \Spin(r,s)\to 
\SO(r,s)$ is the double covering of $\SO(r,s)$ by $\Spin(r,s)$, we have to check, for which of 
the holonomy groups of Berger's list $V_{\Hu}$ is non-zero. 

We obtain
\vab

{\bf Theorem.} \hab {\it Let $M^{r,s}$ be a simply connected, non-locally 
symmetric irreducible semi-Riemannian spin manifold of dimension $n=r+s$ and index $r$. 
Let $\fP$ be the space of parallel spinors of $M$. Then $\dim \fP=N>0$ if and only if the 
holonomy representation $H$ of $M$ is (up to conjugacy in the full orthogonal group) one in 
the list given below. Furthermore, this list gives the chiral type and in the case of 
space- and time oriented spin manifolds the causal type of the parallel spinors. 
\benu 
\item 
$H=\SU(p,q)\subset \SO(2p,2q)$, $n=2(p+q)$, $r=2p$.\\ 
Then $N=2$. There exists a basis $\{\gf_1,\gf_2\}$ of $\fP$ such that\\ 
\mbox{}\hspace{3ex} 
\parbox{12cm}{$\gf_1,\,\gf_2\in \gG(S^+)$ or $\gf_1,\,\gf_2\in \gG(S^-)$ if 
$p+q$ is even,\\ 
$\gf_1\in \gG(S^+)$, $\gf_2\in \gG(S^-)$ if $p+q$ is odd,\\
$\gf_1,\,\gf_2$ have the same causality if $p$ is even\\
$\gf_1,\,\gf_2$ have different causality if $p$ is odd.
}
\item 
$H=\Sp(p,q)\subset \SO(4p,4q)$, $n=4(p+q)$, $r=4p$.\\ 
Then $N=p+q+1$, $\fP\subset\gG(S^+)$ or $\fP\subset \gG(S^-)$. All 
non-trivial parallel spinors have the same causal type. 
\item 
$H=\G_2 \subset \SO(7)$, $n=7$, $r=0$.\\
Then $N=1$.
\item 
$H=\G^*_{2(2)}\subset \SO(4,3)$, $n=7$, $r=4$.\\
Then $N=1$. All non-trivial parallel spinors are non-isotropic.
\item 
$H=\G^{\CC}_2\subset \SO(7,7)$, $n=14$, $r=7$.\\
Then $N=2$. There exists a basis $\{\gf_1,\,\gf_2\}$ of $\fP$ such that 
$\gf_1\in\gG(S^+)$, $\gf_2\in \gG(S^-)$, $\gf_1,\,\gf_2$ are isotropic but not orthogonal to each 
other. 
\item 
$H=\Spin(7)\subset \SO(8)$, $n=8$, $r=0$.\\
Then $N=1$.
\item 
$H=\Spin(4,3)\subset \SO(4,4)$, $n=8$, $r=4$.\\
Then $N=1$. All non-trivial parallel spinors are non-isotropic.
\item 
$H=\Spin(7)^{\CC}\subset \SO(8,8)$, $n=16$, $r=8$.\\
Then $N=1$. All non-trivial parallel spinors are non-isotropic.
\eenu
}

\vab

This theorem shows that 
the possible holonomy group of a manifold with parallel spinors is uniquely 
determined (up to conjugacy) by the dimension and the index of $M$ and the number of linearly 
independent parallel spinor fields.\\ 
In the proof of the theorem we use the fact that the vector spaces $V_{\tilde{H_1}}$ 
and $V_{\tilde{H_2}}$ are isomorphic if the Lie algebras $\tilde{\fh_1}$ and 
$\tilde{\fh_2}$ of $\tilde{H_1}$ and $\tilde{H_2}$, respectively, are different real forms of 
the same complex Lie algebra. Therefore, the dimension of $V_{\tilde{H}}$ is given by 
that of the space corresponding to the compact real form. The explicit description of 
$V_{\tilde{H}}$ allows to determine the causal type of the parallel spinors in the 
pseudo-Riemannian situation . Another proof of the theorem by a straight forward direct calculation 
of the vector space $V_{\tilde{H}}$ using standard formulas of the spinor calculus 
is given in \cite{BK 97}. 

\mk

Section~\ref{s-spinor representations} of this paper contains some necessary 
notations from the spinor calculus. In Section~\ref{s-some facts} we recall 
the complete Berger list
of the holonomy representations of irreducible simply connected non-locally 
symmetric semi-Riemannian manifolds, explain the relation between parallel spinors and 
these holonomy groups and fix some notation concerning the appearing groups. In 
Section~\ref{s-fixed spinors} we calculate the spaces $\VH$ for the holonomy groups of Berger's list 
explicitly
and determine the chiral and the causal type of the corresponding parallel 
spinors. 
\mk

The authors would like to thank R.~Bryant for helpfull remarks concerning the 
original Berger list. 

\section{Spinor representations}\label{s-spinor representations}

For a more detailed explanation of the following facts see e.g. \cite{LM 89}, 
\cite{Ba 81}, \cite{F 97}. 

Let $\cC_{r,s}$ denote the Clifford algebra of $\,\RR^{r,s}= 
(\RR^{r+s},\,\langle\,,\,\rangle_{r,s})\,$, where $\langle\,,\,\rangle_{r,s}$ is a scalar product of index $r$. It is 
generated by an orthonormal basis $e_1,\ldots,e_{r+s}$ of $\RR^{r,s}$ with relations 
$\,e_i\cdot e_j+e_j\cdot e_i=-2\gk_j\cdot\gd_{ij}\,$,
where
$\,\gk_j=\langle e_j\,,\,e_j\rangle_{r,s}=\pm 1\,$.\\
If $n=r+s$ is even, then the complexification $\cC^{{\mathbb C}}_{r,s}$ of 
$\cC_{r,s}$ is isomorphic to the matrix algebra $\CC(2^{\frac n2})$. If $n=r+s$ is odd, then 
$\cC^{{\mathbb C}}_{r,s}$ is isomorphic to $\CC(2^{[\frac n2]}) \oplus \CC(2^{[\frac n2]})$. 

We will use the following isomorphisms:\\
Let $
E=\mbox{{\footnotesize $\Big(\barr{cc}1&0 \\ 0&1\earr\Big)$}}, \; 
T=\mbox{{\footnotesize $\Big(\barr{cc} 0&-i\\ i&0\earr\Big)$}}, \;
U=\mbox{{\footnotesize $\Big(\barr{cc}i&0\\ 0&-i\earr\Big)$}}, \;
V=\mbox{{\footnotesize $\Big(\barr{cc}0&i\\ i&0\earr\Big)$}}$
and
\[
\gt_k=\left\{
\barr{ccl}
i & {\rm if} & \gk_k=-1,\\
1 & {\rm if} & \gk_k=1.
\earr
\right.
\]
In case $n=r+s$ even we define
\[
\gF_{r,s}\,:\,\cC^{{\mathbb C}}_{r,s} \stackrel{\sim}{\longrightarrow} 
\CC(2^{\frac n2}) 
\]
by 
\be\label{phi1}
\barr{r@{\; = \; }l}
\gF_{r,s}(e_{2k-1}) & \gt_{2k-1}\,E\otimes\cdots\otimes E\otimes U\otimes 
\underbrace{T\otimes\cdots\otimes T}_{k-times}\\
\gF_{r,s}(e_{2k}) & \gt_{2k}\,E\otimes\cdots\otimes E\otimes V\otimes 
\underbrace{T\otimes\cdots\otimes T}_{k-times}\,,
\earr
\ee
$k=1,\ldots,\frac n2$.
The restriction of $\gF_{r,s}$ to $\Spin(r,s)\subset \cC_{r,s}$ yields the 
spinor representation $\gD_{r,s}$ which we will use here. 

In case $n=r+s$ odd
\[
\gF_{r,s}\,:\, \cC^{{\mathbb C}}_{r,s} \stackrel{\sim}{\longrightarrow} 
\CC(2^{[\frac n2]})\oplus \CC(2^{[\frac n2]}) 
\]
is given by
\be\label{phi2}
\barr{r@{\; = \; }l}
\gF_{r,s}(e_k) & (\gF_{r,s-1}(e_k)\,,\, \gF_{r,s-1}(e_k)), \;\;\; 
k=1,\ldots, n-1\\[2ex] 
\gF_{r,s}(e_n) & \gt_n (i\,\, T\otimes\cdots \otimes T\,, 
\,-i\,\,T\otimes\cdots \otimes T). 
\earr
\ee
If we restrict $\gF_{r,s}$ to $\Spin(r,s)\subset\cC_{r,s}$ and project onto the 
first component we obtain the spinor representation $\gD_{r,s}$ in the 
odd-dimensional case. 

The Lie algebra $\fspin(r,s)$ of $\Spin(r,s)$ is given by 
\[
\fspin(r,s)= {\rm span}(e_k\cdot e_l \,|\, 1\le k<l\le n)\subset \cC_{r,s}. 
\]
Let $D_{kl}$ be the $(n \times n)$-matrix whose $(k,l)$-entry is 1 and all of whose 
others are 0 and let $E_{kl}$ be the matrix 
\be\label{E}
E_{kl}=-\gk_l D_{kl}+\gk_k D_{lk}\,.
\ee 
The differential of the double covering $\gl\,:\,\Spin(r,s)\to\SO(r,s)$ is 
given by 
\[
\gl_*(e_k\cdot e_l)=2 E_{kl}\,.
\]
Let $u(\gep)\in \CC^2$, $\gep=\pm 1$, be the vectors $u(\gep)=\frac 1{\sqrt{2}} 
\mbox{{\footnotesize $\Big(
\barr{c}
1\\
 -\gep \cdot i
\earr\Big)$}}$ of $\CC^2$. For calculations we will use the orthonormal basis 
\[
\{u(\gep_m,\ldots, \gep_1):= u(\gep_m) \otimes u(\gep_{m-1}) \otimes \cdots 
\otimes u(\gep_1) \; |\; \gep_j=\pm 1\} 
\]
of $\gD_{r,s}=\CC^2 \otimes \cdots \otimes \CC^2 \;\; (m=[\frac n2])$.

In case of even $n=r+s$ the splitting $\gD_{r,s}=\gD^+_{r,s} \oplus 
\gD^-_{r,s}$ of $\gD_{r,s}$ in positive and negative spinors is given by 
\[
\gD^{\pm}_{r,s}=\left\{u(\gep_m, \ldots, \gep_1)\, | \, \prod^m_{j=1} 
\gep_j=\pm 1\right\}. 
\]  

If we are in the pseudo-Euclidean case $(0<r<r+s)$, there exists an indefinite 
scalar product $\langle\,, \rangle$ on $\gD_{r,s}$ which is invariant under the action 
of the connected component $\Spin_0(r,s)$ of $\Spin(r,s)$. $\langle\,,\rangle$ 
is given by 
\[
\langle u\,,\, v \rangle = i^{\frac{r(r-1)}2} (e_{i_1}\cdots e_{i_r}€\cdot 
u\,,\, v)_{{\mathbb C}^{2^m}}, 
\]
where $e_{i_1}, \ldots, e_{i_r}\in \RR^{r,s}$ are the timelike vectors of the 
orthonormal basis \lb $e_1, \ldots, e_n$ of $\RR^{r,s}\,$, $i_1<\cdots <i_r$. 

In particular, if $r$ and $s$ are odd, then $\gD^+_{r,s}$ and $\gD^-_{r,s}$ are 
totally isotropic. If $r$ and $s$ are even, $\gD^+_{r,s}$ is $\langle\,, 
\rangle$-orthogonal to $\gD^-_{r,s}$. 

In the case that
\[
\langle x\,,\,y\rangle_{r,s} = -x_1y_1-\cdots - 
x_ry_r+x_{r+1}y_{r+1}+\cdots +x_{r+s}y_{r+s} 
\]
and that $e_1,\ldots, e_{r+s}$ is the standard basis of $\RR^{r+s}$ we have for 
even index $r=2\hat{r}$ 
\bea\label{sk1}
\lefteqn{\langle u(\gep_m,\ldots, 
\gep_1)\,,\,u(\gd_m,\ldots,\gd_1)\rangle\nonumber}\\ 
&=&\left\{
\barr{cll}
0 & \mbox{if} & (\gep_m,\ldots, \gep_1)\not= (\gd_m,\ldots,\gd_1)\\[1ex] 
\gep_1\cdots \gep_{\hat{r}} & \mbox{if} & (\gep_m,\ldots, \gep_1)= 
(\gd_m,\ldots,\gd_1). 
\earr
\right.
\eea
If we consider the bilinear form
\[
\langle x\,,\,y\rangle_{r,r} = -\sum^r_{j=1} x_{2j-1}\cdot 
y_{2j-1}+\sum^r_{j=1} x_{2j}\cdot y_{2j} 
\]
and the standard basis  $(e_1,\ldots, e_{2r})$ of $\RR^{2r}$ we obtain
\bea\label{sk2}
\lefteqn{\langle u(\gep_r,\ldots, \gep_1)\,,\,u(\gd_r,\ldots 
,\gd_1)\rangle}\nonumber\\ 
&=&\!\!\left\{
\barr{ll}
0 & \mbox{if} \; (\gep_r,\ldots, \gep_1)\not= (-\gd_r,\ldots,-\gd_1)\\[1ex] 
(-i)^{\hat{r}}\, \gep_1\cdot \gep_3\cdots \gep_{r-1} & \mbox{if} \; 
r=2\hat{r} \;\mbox{and}\\ 
& (\gep_r,\ldots, \gep_1)= (-\gd_r,\ldots,-\gd_1)\\[1ex]
-i^{\hat{r}}\, \gep_2\cdot \gep_4\cdots \gep_{r-1} & \mbox{if}\; 
r=2\hat{r}+1 \;\mbox{and}\\ &(\gep_r,\ldots, \gep_1)= (-\gd_r,\ldots,-\gd_1). 
\earr
\right.
\eea
Let $\, j: \cC_{r,s} \hookrightarrow \cC_n^{{\mathbb C}}\,$ be the embedding 
which sends $e_k$ to $e_k$ if $\gk_k=1$ and $e_k$ to $ie_k$ if $\gk=-1$. Then 
$j(\cC_{r,s})$ and $j(\fspin (r,s))$ are real forms of $\cC_n^{{\mathbb C}}$ and $\fspin 
(n)^{{\mathbb C}}$, respectively. We have 
\bea \label{*1} 
\gF_{r,s} = \gF_n \circ j.  
\eea

\section{Some facts on holonomy groups of semi-\\Riemannian spin 
manifolds}\label{s-some facts} 
A connected semi-Riemannian manifold $(M^{r,s},g)$ of dimension $r+s$ and index 
$r$ is called {\it irreducible}, if its holonomy representation is irreducible. It is 
called {\it indecomposable} if its holonomy groups $H_u\subset \O(r,s)$ do not leave 
invariant any nondegenerate proper subspace of $\RR^{r,s}$. Of course, in the Riemannian 
case, indecomposable manifolds are irreducible. In the pseudo-Riemannian case there 
exist holonomy groups whose holonomy representations are not irreducible but which are 
not decomposable into a direct sum of pseudo-Riemannian holonomy representations 
(\cite{Wu 67}, \cite{BI 93}). 

De Rham's splitting theorem reduces the study of complete simply connected 
semi-Riemannian manifolds to indecomposable ones. 

\begin{theorem} {\rm (\cite{DR 52}, \cite{Wu 64})} \hab \label{th1}
Let $(M^{r,s},g)$ be a simply connected complete semi-Riemannian manifold. Then 
$(M^{r,s},g)$ is isometric to a product of simply connected complete indecomposable 
semi-Riemannian manifolds. 
\end{theorem}

Up to now, there is no list of all indecomposable restricted holonomy groups in the 
pseudo-Riemannian setting, but the irreducible cases are known. 

The irreducible pseudo-Riemannian symmetric spaces were classified by \lb 
M.~Berger in \cite{Be 57}. The irreducible restricted holonomy groups of non-locally 
symmetric spaces 
are listed in the following (corrected) Berger list.

\begin{theorem} {\rm (\cite{Be 55}, \cite{S 62}, \cite{A 68},
\cite{BG 72}, \cite{Br 96})} 
\label{th2} Let $(M^{r,s},g)$ be a simply connected irreducible non-locally 
symmetric semi-Riemannian manifold of dimension $n=r+s$ and index $r$. Then its holonomy 
representation is (up to conjugacy in $\O(r,s)$) one of the following 
\[
{\renewcommand{\arraystretch}{1.05}\barr{llllll}
\SO_0(p,q) & & & & n=p+q\ge 2, & r=p\\ 
{\arraycolsep0mm  \barr{l}
\U(p,q)\\
\SU(p,q)
\earr}       & {\arraycolsep0mm  \barr{c}
                  \subset \\
                  \subset
                  \earr}       &  {\arraycolsep0mm  \barr{l} 
                                             \SO(2p,2q)\\
                                             \SO(2p,2q)
\earr} & \left. {\arraycolsep0mm \barr{l} 
\mbox{}\\ 
\mbox{} 
\earr}\right\} & n=2(p+q)\ge 4, & r=2p\\ 
{\arraycolsep0mm  \barr{l}
\Sp(p,q)\\
\Sp(p,q)\cdot \Sp(1)
\earr}      &  {\arraycolsep0mm  \barr{c} 
                  \subset \\
                  \subset
                  \earr}       & {\arraycolsep0mm  \barr{l}
                                     \SO(4p,4q)\\
                                     \SO(4p,4q)
\earr} & \left. {\arraycolsep0mm \barr{l} 
\mbox{}\\ 
\mbox{} 
\earr}\right\} & n=4(p+q)\ge 8, & r=4p\\ 
\Sp(p,\RR)\cdot \SL(2,\RR) & \subset & \SO(2p,2p) && n=4p\ge 8, &r=2p\\ 
\Sp(p,\CC)\cdot \SL(2,\CC) & \subset & \SO(4p,4p) && n=8p\ge 16, &r=4p\\ 
\SO(p,\CC) &\subset &\SO(p,p) & & n=2p\ge 4, & r=p\\ 
{\arraycolsep0mm  \barr{l}
\G_2\\
\G_{2(2)}^*
\earr}         & {\arraycolsep0mm  \barr{c}
                    \subset \\
                    \subset
                    \earr}    & {\arraycolsep0mm  \barr{l}
                                  \SO(7)\\
                                  \SO(4,3) 
\earr} & \left. {\arraycolsep0mm \barr{l} 
\mbox{}\\ 
\mbox{} 
\earr}\right\} & n=7 &\\ 
\G_2^{\CC} &\subset &\SO(7,7) && n=14 &\\ 
{\arraycolsep0mm  \barr{l}
\Spin(7)\\
\Spin(4,3)
\earr}            & {\arraycolsep0mm  \barr{c}
                       \subset\\ 
                       \subset
                       \earr}  & {\arraycolsep0mm  \barr{l}
                                    \SO(8)\\
                                    \SO(4,4)
\earr} & \left. {\arraycolsep0mm \barr{l} 
\mbox{}\\ 
\mbox{} 
\earr}\right\} & n=8 &\\ 
\Spin(7)^{{\mathbb C}} & \subset &\SO(8,8) & &n=16 & 
\earr
}
\]
\end{theorem}

It is known for all groups $H$ of the list in Theorem~\ref{th2} 
that there exists a non-symmetric semi-Riemannian manifold with holonomy $H$ 
(see \cite{Br 87}, \cite{Br 96}). 

Now let $(M^{r,s},g)$ be a semi-Riemannian spin manifold with spin structure 
$(Q,f)$ and spinor bundle $S=Q\times_{\Spin(r,s)}\gD_{r,s}$. 

A spinor field $\gf\in \gG(S)$ is called {\it parallel}, if $\nabla^S\gf\equiv 
0$, where $\nabla^S$ is the spinor derivative associated to the Levi-Civita 
connection $\tilde{A}$ on $Q$. 

Let us denote by $\gt^{\tilde{A}}_{\go}\,:\,Q_x \to Q_x$ the parallel transport 
in the spin structure $Q$ with respect to $\tilde{A}$ along a closed curve $\go$ 
starting at $x\in M$. Then for $\tilde{u} \in Q_x$ 
\bean
\tilde{H}_{\tilde{u}} & :=& \{ g\in \Spin(r,s)\;|\; \mbox{There exists a closed 
curve} \: \go\\ 
&& \hspace{17ex} \mbox{starting at} \: x \: \mbox{such that} \: \tu\cdot g = 
\gt^{\tilde{A}}_{\go}(\tilde{u})\} 
\eean 
denotes the holonomy group of $\tilde{A}$ with respect to $\tilde{u}€\in Q_x$. If 
$\tilde{u}_1 \in Q_{x_1}$ is another point in $Q$ and $\gs$ is a curve from $x$ to $x_1$, then 
\[
\tilde{H}_{\tilde{u}_1}= a^{-1}\cdot \tilde{H}_{\tilde{u}} \cdot a,
\]
where $a€\in \Spin(r,s)$ is the element with $u_1=\gt^{\tilde{A}}_{\gs}(u) 
\cdot a$. Hence, all holonomy groups of $\tilde{A}$ are conjugated to each other. 

Let
\bean
Q^{\tilde{A}}(\tilde{u}) & := & \{ \hat{u}\in Q\;|\; \mbox{There exists an } 
\tA-\mbox{horizontal curve}\\ 
&&\hspace{10ex} \mbox{connecting }\tu \mbox{ with } \hat{u}\}
\eean
be the holonomy bundle of $(Q,\tA)$ with respect to $\tu \in Q$. According to the 
reduction theorem, 
$(Q,\tA)$ reduces to the $\tH_{\tu}$-principal bundle $Q^{\tA}(\tu)$. Hence, 
the spinor bundle $S$ of $M^{r,s}$ is given by 
\[
S=Q\times_{\Spin(r,s),\gr}\gD_{r,s}=Q^{\tA}(\tu)\times_{\tH_{\tu},\gr} 
\gD_{r,s}. 
\]
Then, we have a bijection between the space $\fP$ of all parallel spinor fields of 
the connected manifold $M^{r,s}$ and the space $V_{\tH_{\tu}}$ of all fixed spinors 
of $\gD_{r,s}$ with respect to the holonomy group $\tH_{\tu}$ which is given by 
\[
\barr{lccl}
V_{\tH_{\tu}} :=\!\! & \{ v\in \gD_{r,s}\;|\;\gr(\tH_{\tu})v=v\} & \mapsto & 
\fP\\ 
& v & \mapsto & \gf_v \in \gG(S)\\
&&&\gf_v(y) = [\gt^{\tA}_{\go}(\tu),v], \mbox{ where } \, \go\\ 
&&& \mbox{is a curve in } \, M \, \mbox{ from } \, x \, \mbox{ to } \, y. 
\earr
\]  
Now, let us suppose that $M^{r,s}$ is simply connected. Then the holonomy groups 
$\tH_{\tu}$ coincide with the restricted holonomy groups, in particular, they are connected. 
Therefore, the vector space $V_{\tH_{\tu}}$ equals 
\[
V_{\tilde{\fh}_{\tilde{u}}}:=\{v\in \gD_{r,s}\;|\; 
\gr_*(\tfh_{\tu})v=0\}, 
\]
where $\tfh_{\tu}$ is the Lie algebra of $\tH_{\tu}$. If $H_u \subset \SO(r,s)$ is 
the holonomy group of the Levi-Civita connection $A$ on the bundle of positively 
oriented frames of $M^{r,s}$ with respect to a frame $u$ in $x\in M$ and if $\tu \in Q_x$ 
denotes a lift of $u$ into the spin structure $Q$, then $\gl(\tH_{\tu})=H_u$, where 
$\gl\;:\; \Spin(r,s)\to \SO(r,s)$ is the double covering of $\SO(r,s)$ by $\Spin(r,s)$. 
Hence, we identify the Lie algebra $\fh_u$ of $H_u$ with $\tfh_{\tu}$ using 
$\gl^{-1}_*$. 

In Berger's list (Theorem~\ref{th2}) the holonomy groups are described up to 
conjugacy in the full orthogonal group $\O(r,s)$. Since we consider holonomy groups of 
oriented manifolds each of the conjugacy classes of the groups $H$ of Theorem~\ref{th2} 
gives rise to two different conjugacy classes in $\SO(r,s)$, generated by $H$ and by $H' 
=T_1\cdot H \cdot T^{-1}_1$, where $T_1= \diag(-1,1,\ldots,1)$. It depends on the chosen 
orientation of $M$ which conjugacy class appears. For the Lie algebras $\fh'$ and $\fh$ of $H'$ 
and $H$ we have $\fh'=\Ad(T) \fh$. For the Lie algebras $\tfh'$ and $\tfh$ of the 
corresponding holonomy groups $\tH'$ and $\tH$ of the spin structure $(Q,\tA)$ it follows 
\[
\tfh'=\Ad(\tT) \tfh,
\]
where $\gl(\tT)=T$. Since $\tT=\pm e_1\in \fpin(r,s)\subset \cC_{r,s}$ and 
$\Ad(\tT)X=\tT\cdot X\cdot \tT^{-1}\in \cC_{r,s}$ for all $X\in \cC_{r,s}$ it results 
$\tfh'=-\gk_1e_1€\cdot \tfh \cdot e_1 \subset \cC_{r,s}$. Therefore we have 
\[
V_{\tH'}=e_1\cdot V_{\tH}\,.
\]    
Hence, it is sufficient to calculate $V_{\tH}$. By Clifford multiplication of 
$v\in V_{\tH}$ with $e_1$ the chiral type of $v$ is changed. 

If we want to study the causal type of parallel spinors we have to restrict ourselves 
to space- and time-oriented pseudo-Riemannian manifolds in order to have an 
indefinite scalar product on the spinor bundle $S$. 
In that case to each conjugacy class of the group $H$ of Theorem~\ref{th2} 
correspond four different conjugacy classes in $\SO_0(r,s)$ generated by 
\[
H,\; H'=T_1HT^{-1}_1, \; H''=T_2HT^{-1}_2,\; H'''=T_3HT_3^{-1},
\]
where $T_1=\diag(-1,1,\ldots,1)$, $T_2=\diag(1,\ldots,1,-1)$,\\ 
$T_3=\diag(-1,1,\ldots,1,-1)$. With analogous notations as above we obtain 
\bean
V_{\tH'} & = & e_1\cdot V_{\tH}\\
V_{\tH''} & = & e_n \cdot V_{\tH}\\
V_{\tH'''} & = & e_1 \cdot e_n \cdot V_{\tH}.
\eean
(Here $e_1$ is timelike and $e_n$ is spacelike.)
In the first two cases the chiral character of $v\in V_{\tH}$ is changed, in the 
third case it remains the same. Since 
\bean
\langle e_j\cdot v\,,\,e_j\cdot v\rangle & = & (-1)^{r+1}€\langle v\,, \, e_j 
\cdot e_j \cdot v\rangle\\ 
&=& (-1)^r \gk_j\langle v\,,\,v\rangle\qquad j=1,\ldots,n,\, v\in \gD_{r,s}, 
\eean
in case of even (odd) index $r$ the causal type of $v€\in V_{\tH}$ is changed 
(remains the same) by Clifford multiplication with $e_1$ or $e_1\cdot e_n$ and remains the 
same (is changed) by Clifford multiplication with $e_n$. 
\mk

Next we describe the groups which occur in Theorem~\ref{th2} in more detail. We 
identify 
\bean
\HH^n & \cong & \CC^{2n}\\
\left(\barr{c}
a_1\\
\mbox{}\\
\vdots\\
\mbox{}\\
a_n
\earr\right) & \mapsto & \left(\barr{c}
z_1\\
\bar{w}_1\\
\vdots\\
z_n\\
\bar{w}_n
\earr\right), \qquad \mbox{where } \, a_{\gn}=z_{\gn}+w_{\gn}€\cdot j\qquad 
(\gn=1,\ldots,n) 
\eean
and
\bean
\CC^n & \cong & \RR^{2n}\\
\left(\barr{c}
z_1\\
\mbox{}\\
\vdots\\
\mbox{}\\
z_n
\earr\right) & \mapsto & \left(\barr{c}
x_1\\
y_1\\
\vdots\\
x_n\\
y_n
\earr\right), \qquad \mbox{where } \, z_{\gn}=x_{\gn}+iy_{\gn}\qquad 
(\gn=1,\ldots,n). 
\eean
We denote by $J_{{\mathbb R}}$ the matrix
\[
J_{{\mathbb R}}=
\left(\barr{ccccc}
J        & 0   & \cdots &   & 0\\
0        & J   &            &   &  \\
\vdots &     & \ddots  &  & \vdots\\
           &     &            &   & 0\\
0         &      &\cdots & 0 & J 
\earr\right)
\in \GL(2n,\RR),
\]
where $J=\mbox{{\footnotesize $\Big(\barr{cc}
0&  -1\\ 
1& 0\earr\Big)$}}$
and by $J_{{\mathbb C}}$ the same matrix in $\GL(2n,\CC)$. Let $\KK(n)$ be the 
space of $(n\times n)$-matrices with entries in the (skew-) field $\KK$. We consider the 
embeddings $i_{{\mathbb C}}$ and $i_{{\mathbb R}}$ associated to the identifications 
$\HH^n\cong \CC^{2n}$ and $\CC^n \cong \RR^{2n}:$ 
\[
i_{{\mathbb C}}\;:\; \HH(n) \hookrightarrow \CC(2n)
\]
\[
i_{{\mathbb 
C}}\left((a_{\gm\gn})_{\gm,\gn=1,\ldots,n}\right)=\left(\left(\barr{cc} 
z_{\gm\gn} & -w_{\gm\gn}\\
\bar{w}_{\gm\gn} & \bar{z}_{\gm\gn}
\earr\right)\right)_{\gm, \gn=1,\ldots,n},
\] 
where $a_{\gm\gn}=z_{\gm\gn}+w_{\gm\gn}\cdot j$,
\[
i_{{\mathbb R}}\;:\; \CC(n) \hookrightarrow \RR(2n)
\]
\[
i_{{\mathbb R}}\left((z_{\gm\gn})_{\gm,\gn=1,\ldots,n}\right)
=\left(\left(\barr{cc}
x_{\gm\gn} & -y_{\gm\gn}\\
y_{\gm\gn} & x_{\gm\gn}
\earr\right)\right)_{\gm, \gn=1,\ldots,n},
\] 
where $z_{\gm\gn}=x_{\gm\gn}+iy_{\gm\gn}$.

Then we have
\bean
i_{{\mathbb C}}(\HH(n)) &= & \{A\in \CC(2n)\,|\;\bar{A} J_{{\mathbb 
C}}=J_{{\mathbb C}}A\}\\ 
i_{{\mathbb R}}(\CC(n)) &= & \{A\in \RR(2n)\,|\;A J_{{\mathbb R}}=J_{{\mathbb 
R}}A\}. 
\eean
Furthermore, let $I^{{\mathbb K}}_{p,q}$ be the matrix
\[
I^{{\mathbb K}}_{p,q}=\left(\barr{cc}
-E_p & 0\\
0& E_q
\earr
\right)
\in \GL(p+q,\KK),
\]
where $E_r$ denotes the unity matrix in $\KK(r)$.

The special pseudo-orthogonal group
\[
\SO(p,q)=\{A\in \SL(p+q,\RR)\;|\; A^TI^{{\mathbb R}}_{p,q}A=I^{{\mathbb 
R}}_{p,q}\} 
\]
(in its standard form) is the invariance group of the bilinear form
\[
\langle 
x\,,\,y\rangle=-x_1y_1-\cdots-x_py_p+x_{p+1}y_{p+1}+\cdots+x_{p+q}y_{p+q} 
\]
on $\RR^{p+q}$.

We identify the pseudo-unitary group
\[
U(p,q)=\{A \in \CC(p+q)\;|\; \bar{A}^T I^{{\mathbb C}}_{p,q} A=I^{{\mathbb 
C}}_{p,q}\} 
\]
with the subgroup $i_{{\mathbb R}}(U(p,q))=\SO(2p,2q)\cap i_{{\mathbb 
R}}\CC(p+q)$ of $\SO(2p,2q)$ and the symplectic group 
\[
\Sp(p,q)=\{A \in \HH(p+q)\;|\; 
\bar{A}
^T I^{{\mathbb H}}_{p,q} A= I^{{\mathbb H}}_{p,q}\}
\]
with the subgroup $i_{{\mathbb C}}€\Sp(p,q)=U(2p,2q)\cap i_{{\mathbb 
C}}(\HH(p+q))$ of $U(2p,2q)$. 

Now we are going to describe the subgroup $\Sp(p,q)\cdot \Sp(1)$ of $\SO(4p,4q)$. 
Each quaternion $a\in \Sp(1)$ defines an orthogonal map $R_a\,:\,\RR^{4p+4q} \to 
\RR^{4p+4q}$ by right multiplication with $a$ on $\HH^{p+q}=\RR^{4p+4q}$. In particular, if 
$a=x_0+iy_0+(x_1+iy_1)j$ and \[ 
r_a:=\left(
\barr{cccc}
x_0 & -y_0 &-x_1 & y_1\\
y_0 & x_0 & y_1 & x_1\\
x_1 & -y_1 & x_0 & -y_0\\
-y_1 & -x_1 & y_0 & x_0
\earr\right)
\]
then we have 
\[
R_a= \left(
\barr{ccccc}
r_a & 0 & \cdots && 0\\
0 & r_a &&&\\
\vdots &&\ddots & & \vdots\\
&&&& 0\\
0 &&\cdots & 0 & r_a
\earr
\right)\in \SO(4p,4q).
\]
The group $\Sp(p,q)\cdot \Sp(1)$ equals
\[
Sp(p,q)\cdot \Sp(1)=\{A\cdot R_a\;|\; A\in i_{{\mathbb R}} i_{{\mathbb C}} 
\Sp(p,q),\, a\in \Sp(1)\}. 
\]
Let $\go_\KK$ be a non-degenerate skew-symmetric bilinear form on $\KK^{2p}$, 
$\KK=\RR,\, \CC$. By $\Sp(p,\KK)$ we denote the group of all automorphisms of 
$(\KK^{2p},\go_\KK)$. The bilinear form $g$ on $\RR^{4p}=\RR^{2p}\otimes \RR^2$ defined by 
\[
g(x\otimes a, y\otimes b) := \det(a,b)\cdot \go_\RR(x,y)
\]
is a metric of signature $(2p,2p)$ on $\RR^{4p}$. Hence, 
$\Sp(p,\RR)\cdot\SL(2,\RR)$ can be considered as a subgroup of all automorphisms of $(\RR^{4p},g)$, where 
$A\cdot B\in \Sp(p,\RR)\cdot \SL(2,\RR)$ acts on $\RR^{4p}=\RR^{2p}\otimes \RR^2$ by 
\[
A\cdot B(x\otimes a) =Ax\otimes Ba.
\]
In the same way, the bilinear form $h$ on $\RR^{8p}=\CC^{4p}=\CC^{2p}\otimes 
\CC^2$ defined by 
\[
h(z\otimes a,w\otimes b)=\Re(\det(a,b)\cdot \go_\CC(z,w))
\]
gives an embedding of $\Sp(p,\CC)\cdot \SL(2,\CC)$ into $\SO(4p,4p)$. 

The remaining classical group
of the list in Theorem~\ref{th2} is
\[
\SO(n,\CC)  =  \SL(n,\CC)\cap \{A\in \CC(n)\;|\;A^TA=E_n\}.
\]
Let us come to the exceptional cases. $\G_2$, $\G^*_{2(2)}$, $\Spin(7)$, 
$\Spin(4,3)$ are subgroups of $\SO(7)$, $\SO(4,3)$, $\SO(8)$ and $\SO(4,4)$, respectively. 
They are defined as invariance groups of so-called €{\it nice $3$-forms} and 4-{\it 
forms}, respectively. Let $w_0,\, w_1\in \gL^3(\RR^7)$, $\gs_0,\,\gs_1 \in 
\gL^4(\RR^8)$ be the forms 
\bean
w_0 & = & w^{127}+w^{135}-w^{146}-w^{236}-w^{245}+w^{347}+w^{567}\\
w_1 & = & -w^{127}-w^{135}+w^{146}+w^{236}+w^{245}-w^{347}+w^{567}\\
\gs_0 & = & 
\gs^{1234}+\gs^{1256}-\gs^{1278}+\gs^{1357}+\gs^{1368}+\gs^{1458}-\gs^{1467}\\ 
&& 
-\gs^{2358}+\gs^{2367}+\gs^{2457}+\gs^{2468}-\gs^{3456}+\gs^{3478}+\gs^{5678}\\ 
\gs_1 & = & 
\gs^{1234}-\gs^{1256}+\gs^{1278}-\gs^{1357}-\gs^{1368}-\gs^{1458}+\gs^{1467}\\ 
&& 
+\gs^{2358}-\gs^{2367}-\gs^{2457}-\gs^{2468}+\gs^{3456}-\gs^{3478}+\gs^{5678}, 
\eean
where $w^{\ga\gb\gg}=w^{\ga}€\wedge w^{\gb}€\wedge w^{\gg}$ and 
$\gs^{\ga\gb\gg\gd}=\gs^{\ga} \wedge\gs^{\gb} \wedge\gs^{\gg} \wedge\gs^{\gd}$ with respect to the dual 
bases $w^1,\ldots, w^7$ of $e_1,\ldots,e_7\in \RR^7$ and $\gs^1,\ldots, \gs^8$ of 
$e_1,\ldots, e_8\in \RR^8$. 

Then
{\arraycolsep0mm
\bean
&&\G_2 =  \{A\in \SO(7)\;|\;A^*w_0=w_0\}\\
&&\G_{2(2)}^*  =  \{A\in \SO(4,3)\;|\; A^*w_1=w_1\}\\
&&\Spin(7)=\{ A\in \SO(8)\;|\;A^*\gs_0=\gs_0\}\\
&&\Spin(4,3)=\{A \in \SO(4,4)\;|\;A^*\gs_1=\gs_1\}
\eean}\mbox{}\hspace{-5pt}
(compare \cite{Br 87}).

Finally, the groups $\G_2^{{\mathbb C}}\subset \SO(7,\CC)$ and 
$\Spin(7,\CC)\subset\SO(8,\CC)$ are the complexifications of $\G_2\subset \SO(7)$ and $\Spin(7)\subset 
\SO(8)$, respectively. 

\section{The fixed spinors of the holonomy representation}\label{s-fixed 
spinors} 

Let $M^{r,s}$ be a simply connected semi-Riemannian spin manifold of index $r$ and 
dimension $r+s=n$ with holonomy representation $H$ and let $\fh$ be the Lie algebra of $H$. 
Then the parallel spinors are given by the kernel of the representation of 
$\tfh:=\gl^{-1}_*(\fh)$ on the spinor module $\gD_{r,s}$. Hence, we have to check the groups in the 
Berger-Simons list and to determine the kernel 
\[
V_{\tfh}=\{v \in \gD_{r,s}\;|\; \tX \cdot v = 0 \, \mbox{ for any }\, \tX \in \tfh\}. 
\]
We make use of the following obvious fact. If $\fa$ is a subalgebra of $\fspin(n)$ 
and $\fa'$ a subalgebra of $\fspin(r,s)$ such that $j(\fa')^{{\mathbb C}} = 
\fa^{{\mathbb C}} \subset \fspin(n)^{{\mathbb C}}$ then we have by (\ref{*1}) and Weyl's 
unitary trick $V_{\fa}= V_{\fa'}$. 

\subsection{$H=\SO_0(p,q) \;\; (p=r,\,q=s)$}

The spinor representation is either irreducible or decomposes into two 
irreducible representations of the same dimension. Therefore, $V_{\fspin(r,s)} = \{0\}$. 
Consequently, there is no parallel spinor on $M^{p,q}$.

\subsection{$H=\U(p,q)\;\;(r=2p,\, s=2q)$}

We consider $\fh=i_{{\mathbb R}^*}\fu(p,q) \subset \fso(2p,2q)$ and 
$\,\fh_0=i_{{\mathbb R}^*}\fu(\frac{n}{2}) \subset \fso(n)$. Then $\,j(\tilde{\fh})^{{\mathbb 
C}} = \tilde{\fh}_0^{{\mathbb C}}\,$ holds. It is known that $\, 
V_{\fu(\frac{n}{2})} = \{0\} $ (see \cite{Wa 89}). Hence, $V_{\fu(p,q)}=\{0\}$. 

\subsection{$H=\SU(p,q)\;\; (r=2p,\,s=2q)$}

Now we have $\,\fh=i_{\mathbb R^*}\fsu(p,q)\subset \fso(2p,2q)\,$ and 
$\,\fh_0=i_{\mathbb R^*}\fsu(\frac{n}{2}) \subset \fso (n)\,$. Again, $j(\tfh)^{\mathbb C}= 
\tfh_0^{\mathbb C}\,$ holds. By \cite{Wa 89} the dimension of $V_{\fh_0}$ equals 2, thus, 
dim$V_{\fsu(p,q)}=2$. The Lie subalgebra 
$i\,_{\mathbb R^*}\fsu(\frac{n}{2}) \subset \fso(n)\,$ is spanned by $X_{kl}, 
Y_{kl}$ and $D_1-D_k \,\,\, (2 \leq k \leq \frac{n}{2})$, where 
\[
\barr{lcll}
X_{kl}  &:= &  E_{2k-1\,2l-1}+E_{2k\,2l} & (1\le k<l\le \frac{n}{2})\\
Y_{kl} &:= & E_{2k-1\,2l}-E_{2k\, 2l-1} & (1\le k< l\le \frac{n}{2})\\
D_k &:= & E_{2k-1\,2k} & (1\le k\le \frac{n}{2}).
\earr
\]
Using this and (\ref{phi1}) we obtain
\[
V_{\fsu(p,q)} = \mbox{span} \{ u(1,1, \dots ,1), u(-1,-1, \dots , -1) \} 
\]
since both generators are annihilated by $\fsu(\frac{n}{2})$.\\
If $p+q=\frac12 \dim M$ is even then $V_{\fsu(p,q)}\subset \gD^+_{2p,2q}$. If 
$p+q=\frac12 \dim M$ is odd, we have an 1-dimensional space of parallel spinors in $S^+$ as well as 
in $S^-$. According to (\ref{sk1}) the quadratic length of $u(\gep,\ldots,\gep)$ 
is 
\[
\langle u(\gep,\ldots,\gep)\,,\,u(\gep,\ldots,\gep)\rangle=\gep^p.
\]
Hence $u(1,\ldots,1)$ and $u(-1,\ldots,-1)$ are of the same causal type 
(spacelike) if $p$ is even and of opposite causal type (timelike) if $p$ is odd. 

\subsection{$H=\Sp(p,q)\;\;(r=4p,\, s=4q)$}

In this case $\fh$ equals $\,i_{\RR_*}i_{\CC_*}\fsp(p,q) \subset 
\fso(4p,4q)\,$ and 
$\fh_0$ is $\,i_{\RR_*}i_{\CC_*}\fsp(\frac{n}{4})\subset \fso(n)\,$. Then 
$j(\tfh)^{{\mathbb C}}= \tfh_0^{{\mathbb C}}$. According to \cite{Wa 89} the dimension of 
$V_{\fh_0}$ equals $\frac{n}{4}+1$. Hence, dim$V_{\fsp(p,q)}= p+q+1$.\\ 
The subalgebra $\,\fh=i_{\RR_*}i_{\CC_*}\fsp(\frac{n}{4})\subset 
\fso(n)\,$ is spanned by 
{\arraycolsep0mm
\bean
&&X_{2k-1\,2l-1}+X_{2k\,2l},\; Y_{2k-1\,2l-1}-Y_{2k\,2l},\\
&&X_{2k-1\,2l}-X_{2k\,2l-1},\; Y_{2k-1\,2l}+Y_{2k\,2l-1},\;\; (1\le k< l 
\le {\textstyle \frac{n}{4}})\\ 
&&X_{2k-1\,2k},\; Y_{2k-1\,2k}\;\; (1\le k\le {\textstyle \frac{n}{4}})\\ 
&&D_{2k-1}-D_{2k}\;\; (1\le k\le {\textstyle \frac{n}{4}}).
\eean}\mbox{}\hspace{-5pt}
Using this and (\ref{phi1}) one proves that the spinors
\[
\gf_k=\sum_{
\barr{c}
\scst{\gep_1,\ldots,\gep_{p+q}}\\[-0.5ex]
\scst{\gep_i=-1 \: \mbox{{\scsi exactly}}\: k\:\mbox{{\scsi times}}}
\earr}
u(\gep_{p+q},\gep_{p+q},\ldots,\gep_1,\gep_1) 
\]
k=0,1, \dots ,p+q, are annihilated by $\fsp(\frac{n}{4})$. Consequently, they 
constitute a basis of $V_{\fsp(p,q)}$. Obviously $V_{\fsp(p,q)}\subset \gD^+$. 
Furthermore all $\gf_k\; (k=0,\ldots, p+q)$ have the same causal type since $\langle 
\gf_k,\gf_k\rangle= \mbox{{\footnotesize $\Big(\barr{cc} 
p+q\\
k
\earr\Big)$}}>0$. 

\subsection{$H=\Sp(p,q)\cdot \Sp(1)\;\, (r=4p,\,s=4q), \\ H= Sp(p,\RR) \cdot 
SL(2,\RR) \;\, (r=2p,s=2p) $} 

The Lie algebras $\,\fh_1=i_{\RR^*} i_{\CC^*} (\fsp(p,q)\oplus \fsp(1)) 
\subset \fso(4p,4q)\,$, $\fh_2=i_{\RR^*} i_{\CC^*} (\fsp(\frac{n}{4})\oplus 
\fsp(1)) \subset \fso(n)\,$ and $\,\fh_3= \fsp(p,\RR) \oplus \fsl(2,\RR) \subset 
\fso(2p,2p)\,$ are real forms of the complex Lie algebra $\fsp(p,\CC) \oplus \fsl(2,\CC)$. 
Therefore, the vector spaces $V_{\fh_1}, V_{\fh_2}$ and $V_{\fh_3}$ are isomorphic. From 
\cite{Wa 89} we know that $V_{\fsp(\frac{n}{4})\oplus \fsp(1)} = \{0\}$. Hence, there are 
no parallel spinors for $H=Sp(p,q)\cdot Sp(1)$ and $H=Sp(p,\RR) \cdot SL(2,\RR)$. 

\subsection{$H=\G_2,\;G_{2(2)}^*$}

Here $\fh=\fg_{2(2)} \subset \fso(4,3), \,\,\fh_0=\fg_2 \subset \fso(7),\,\, 
j(\tfh)^ {{\mathbb C}}= \tfh_0^{{\mathbb C}}\,$. It is well known that dim$V_{\fg_2}=1$ 
and, therefore, dim$V_{\fg_{2(2)}}=1$. The Lie algebra $\fg_2 
 \subset \fso(7)$ is spanned by 
{\arraycolsep0mm
\bean
&&E_{12}-E_{34},\; E_{12}-E_{56},\; E_{13}+E_{24},\; E_{13}-E_{67},\; 
E_{14}-E_{23},\\ 
&& E_{14}-E_{57},\; E_{15}+ E_{26},\;E_{15}+E_{47},\; E_{16}-E_{25},\; 
E_{16}+E_{37},\\ 
&& E_{17}-E_{36},\;E_{17}-E_{45},\; E_{27}-E_{35},\; E_{27}+E_{46}.
\eean}\mbox{}\hspace{-5pt}
Using this and (\ref{phi2}) one checks that the spinor $\,\varphi= u(1,1,1)+i 
u(-1,-1,-1)\,$ is annihilated by $\fg_2$. Consequently, 
\[ V_{\fg_2} = \;\mbox{span}\{\varphi\} \subset \gD_7\, \hspace{1cm} 
V_{\fg_{2(2)}} =\;\mbox{span}\{\varphi\} \subset \gD_{4,3}. 
\] 
From $\,\langle \varphi,\varphi \rangle =2\,$ we see, that all non-zero spinors 
in $V_{\fg_{2(2)}}$ are spacelike. 

\subsection{$H=\Spin(7)$, $\Spin(4,3)$}

Now, $\fh=\fspin(4,3) \subset \fso(4,4), \,\, \fh_0=\fspin(7) \subset 
\fso(8), \,\, j(\tfh)^{{\mathbb C}} = \tfh_0^{{\mathbb C}}\,$. It is well known that 
dim$V_{\fspin(7)}=1$. Consequently, dim$V_{\fspin(4,3)}=1$. The Lie subalgebra $\fspin(7) \subset 
\fso(8)$ is generated by 
{\arraycolsep0mm
\bean
&&E_{12}+E_{34},\; E_{13}-E_{24},\; E_{14}+E_{23}, \; E_{56}+E_{78},\; 
-E_{57}+E_{68},\; E_{58}+E_{67},\\ 
&& -E_{15}+ E_{26},\; E_{12}+E_{56},\; E_{16}+E_{25}, \; E_{37}-E_{48},\; 
E_{34}+E_{78},\;E_{38}+E_{47},\\ 
&& E_{12}-E_{78},\; E_{17}+E_{28},\; E_{18}-E_{27},\; E_{34}-E_{56}, 
\;E_{35}+E_{46},\; E_{36}-E_{45},\\ 
&& E_{18}+E_{36},\; E_{17}+E_{35},\; E_{26}-E_{48},\; E_{25}+E_{38},\; 
E_{23}+E_{67},\; E_{24}+E_{57}. 
\eean}\mbox{}\hspace{-5pt}
Hence by (\ref{phi1}) the spinor $\,\psi = u(1,-1,-1,1) - u(-1,1,1,-1)\,$ is 
annihilated by $\fspin(7)$. This shows 
\[
V_{\fspin(7)} = \,\mbox{span}\{\psi\} \subset \gD_8^+ , \hspace{1cm} 
V_{\fspin(4,3)} = \,\mbox{span}\{\psi\} \subset \gD_{4,4}^+. 
\]
Furthermore, in case of $\fspin(4,3)$ the spinor $\psi$ is timelike.


\subsection{$H=\SO(p,\CC)\;\;\;(r=p,\,s=p)$}

Now let (in difference to the previous sections) $\SO(p,p)\subset \GL(2p,\RR)$ 
be the invariance group of the inner product $\langle x\,,\,y 
\rangle_{p,p}=-x_1y_1+x_2y_2-x_3y_3+\cdots +x_{2p}y_{2p}$. According to that we have 
\[
\gk_k=\left\{
\barr{cl}
1& k\;\mbox{even}\\
-1& k\; \mbox{odd}
\earr 
\right.
\hspace{5ex}  
\gt_k=\left\{
\barr{cl}
1& k\;\mbox{even}\\
i& k\; \mbox{odd}.
\earr 
\right.  
\]
Then $i_{\RR}(\SO(p,\CC))$ is contained in $\SO(p,p)$. The Lie algebra 
$\,\fso(p) \oplus \fso(p) \subset \fso(2p)\,$ is the compact real form of 
$\,(i_{\RR^*}\fso(p,\CC))^{{\mathbb C}}\,$. Hence, we obtain $\,\,\mbox{dim}\,V_{\fso(p,{\mathbb C})} =$ dim$ 
V_{\fso(p)\oplus \fso(p)}\,$.\\ 
Let $p=2k$. Then the restriction of the $\fso(4k)$-representation 
$\gD_{4k}^\pm$ to $\fso(2k) \oplus \fso(2k)$ is equivalent to the sum $\,\gD_{2k}^\pm \otimes 
\gD_{2k}^+ 
\oplus \gD_{2k}^\mp \otimes \gD_{2k}^-\,$ of tensor products of the 
$\fso(2k)$-representations $\gD_{2k}^\pm$. Since $\,V_{\fso(2k)}= \{0\}\,$ we obtain that 
$\,V_{\fso(2k,{\mathbb C})} = \{0\}\,$.\\ 
In case $p=2k+1$ the restriction $\,\gD_{4k+2}^\pm |_{\fso(2k+1) \oplus 
\fso(2k+1)}\,$ is isomorphic to the tensor product $\,\gD_{2k+1} \otimes \gD_{2k+1}\,$ of the 
$\fso(2k+1)$-representation $\gD_{2k+1}$ by itself. Analogously, we deduce from 
$\,V_{\fso(2k+1)}=\{0\}\,$ that $\,V_{\fso(4k+2,{\mathbb C})}=\{0\}\,$.\\ 
Hence, there are no parallel spinors for $H=SO(p,{\mathbb C})$. 

\subsection{$H=\G^{\CC}_2 \;\;\; (r=7,\,s=7)$}

The compact real form of $\,(i_{\RR^*}(\fg_2^{{\mathbb C}}))^{{\mathbb C}}\,$ 
is equal to $\,\fg_2 \oplus \fg_2 \subset \fso(7) \oplus \fso(7)\,$. Hence, 
$\,\,\mbox{dim} V_{\fg_2^{{\mathbb C}}} = \,\mbox{dim} V_{\fg_2 \oplus\fg_2}\,$. The 
representations $\,\gD_{14}^\pm|{\fg_2 \oplus \fg_2}\,$ of $\,\fg_2 \oplus\fg_2\,$ are 
equivalent to $\,\gD_7|_{\fg_2} \otimes \gD_7|_{\fg_2}\,$. Since $\,\,\mbox{dim} 
V_{\fg_2}=1\,$ we have $\,\,\mbox{dim} V_{\fg_2 \oplus \fg_2}=2\,$, where one parallel spinor 
lies in $\gD_{14}^+$ and the other one in $\gD_{14}^-$. Consequently, 
$\,\,\mbox{dim}\, V_{\fg_2^{{\mathbb C}}}= 2$.\\ 
The Lie algebra   
$i_{\RR^*}\fg^{\CC}_2 \subset i_{\RR^*}\fso(7,\CC) \subset \fso(7,7)$ is 
spanned by 
{\arraycolsep0mm
\bean
&& \xi_{12}-\xi_{34},\; \xi_{13}+\xi_{24},\; \xi_{14}-\xi_{23},\; 
\xi_{12}-\xi_{56},\; \xi_{16}-\xi_{25},\; \xi_{15}+\xi_{26},\; \xi_{36}-\xi_{45},\\ 
&& \xi_{35}+\xi_{46}, \; -\xi_{13}+\xi_{67}, \;-\xi_{14}+\xi_{57},\; 
\xi_{16}+\xi_{37},\; \xi_{15}+\xi_{47},\; \xi_{36}-\xi_{17},\; \xi_{35}- \xi_{27},\\ 
&& \eta_{12}-\eta_{34},\; \eta_{13}+\eta_{24},\; \eta_{14}-\eta_{23},\; 
\eta_{12}-\eta_{56}, \; \eta_{16}-\eta_{25},\; \eta_{15}+\eta_{26},\; \eta_{36}-\eta_{45},\\ 
&& \eta_{35}+\eta_{46},\; - \eta_{13}+\eta_{67},\; -\eta_{14}+\eta_{57},\; 
\eta_{16}+\eta_{37},\; \eta_{15}+\eta_{47},\; \eta_{36}-\eta_{17},\; \eta_{35}-\eta_{27}. 
\eean}\mbox{}\hspace{-5pt}
Then a direct calculation shows that 
\bean
\psi_1 &=& u(1,1,1,1,1,1,1)+u(1,1,1,-1,-1,-1,-1)\\
&& + u(1,-1,-1,-1,-1,1,1)+u(1,-1,-1,1,1,-1,-1)\\
&& + u(-1,-1,1,-1,1,-1,1)-u(-1,1,-1,1,-1,-1,1)\\
&& - u(-1,1,-1,-1,1,1,-1)-u(-1,-1,1,1,-1,1,-1)\\
& \in & \gD_{7,7}^+\\[0.2cm]
\psi_2 &=& u(-1,1,1,1,1,-1,-1)+u(-1,1,1,-1,-1,1,1)\\
&& + u(-1,-1,-1,1,1,1,1)+u(-1,-1,-1,-1,-1,-1,-1)\\
&& -u(1,1,-1,1,-1,1,-1)+u(1,-1,1,-1,1,1,-1)\\
&& + u(1,-1,1,1,-1,-1,1)+u(1,1,-1,-1,1,-1,1)\\
& \in & \gD_{7,7}^-
\eean
are generators of $V_{\fg_2^{{\mathbb C}}}$.\\
From (\ref{sk2}) one calculates $\,\langle \psi_1,\psi_1 \rangle = \langle 
\psi_2,\psi_2 \rangle =0 \,$ and $\,\langle \psi_1,\psi_2 \rangle =8i$. 

\subsection{$H=\Spin(7)^{\CC}\;\;\; (r=8,\,s=8)$}

The Lie algebra $\,(i_{\RR^*}\fspin(7)^{{\mathbb C}})^{{\mathbb C}} \subset 
(i_{\RR^*}\fso(8,\CC))^{{\mathbb C}}\,$ has the compact real form $\,\fspin(7) \oplus \fspin(7) \subset \fso(8) 
\oplus \fso(8) \subset \fso(16)\,$. Hence, $\,\,\mbox{dim}\,V_{\fspin(7)^{\CC}} = 
\,\mbox{dim}V_{\fspin(7)\oplus \fspin(7)}\,$. The representation $\,\gD_{16}^\pm|_{\fspin(7) \oplus \, 
\fspin(7)}\,$ is equivalent to $\,\gD_8^\pm|_{\fspin(7)} \otimes \gD_8^+|_{\fspin(7)} 
\oplus \,\gD_8^\mp|_{\fspin(7)} \otimes\, \gD_8^-|_{\fspin(7)}\,$. Since 
$\,\,\mbox{dim}\,V_{\fspin(7)}=1\,$ and $\,V_{\fspin(7)}\subset \gD_8^+\,$ we obtain 
$\,\,\mbox{dim}\,V_{\fspin(7) \oplus \fspin(7)} =1\,$ and $\,V_{\fspin(7) \oplus \fspin(7)} \subset 
\gD_{16}^+\,$. In particular, $\,\,\mbox{dim}\,V_{\fspin(7)^{{\mathbb C}}} =1\,$. \\ 
The Lie algebra $\,i_{\RR^*}\fspin(7,\CC) \subset \fso(8,8)\,$ is spanned by 
{\arraycolsep0mm
\bean
&& \xi_{12}+\xi_{34},\; \xi_{13}-\xi_{24},\; \xi_{14}+\xi_{23},\; 
\xi_{56}+\xi_{78},\;- \xi_{57}+\xi_{68},\; \xi_{58}+\xi_{67},\; -\xi_{15}+\xi_{26},\\ 
&& \xi_{12}+\xi_{56},\; \xi_{16}+\xi_{25}, \; \xi_{37}-\xi_{48},\; 
\xi_{38}+\xi_{47},\; \xi_{17}+\xi_{28},\; \xi_{18}-\xi_{27},\; \xi_{35}+\xi_{46},\\ 
&& \xi_{36}- \xi_{45},\; \xi_{18}+\xi_{36}, \;\xi_{17}+\xi_{35},\; 
\xi_{26}-\xi_{48},\; \xi_{25}+\xi_{38}, \;\xi_{23}+\xi_{67},\; \xi_{24}+\xi_{57},\\ 
&& \eta_{12}+\eta_{34},\; \eta_{13}-\eta_{24},\; \eta_{14}+\eta_{23},\; 
\eta_{56}+\eta_{78},\; -\eta_{57}+\eta_{68},\; \eta_{58}+\eta_{67},\; -\eta_{15}+\eta_{26},\\ 
&& \eta_{12}+\eta_{56},\; \eta_{16}+\eta_{25},\; \eta_{37}-\eta_{48},\; 
\eta_{38}+\eta_{47},\; \eta_{17}+\eta_{28},\; \eta_{18}-\eta_{27},\; \eta_{35}+\eta_{46},\\ 
&& \eta_{36}-\eta_{45},\; \eta_{18}+\eta_{36},\; \eta_{17}+\eta_{35},\; 
\eta_{26}-\eta_{48},\; \eta_{25}+\eta_{38},\; \eta_{23}+\eta_{67},\; \eta_{24}+\eta_{57}. 
\eean}\mbox{}\hspace{-5pt}
A direct calculation now shows $\,V_{\fspin(7)^{\CC}}\,$ is generated by 
\bean
\eta &=&  u(1,1,1,1,1,1,1,1) -u(1,1,1,1,-1,-1,-1,-1)\\
&&-u(-1,-1,-1,-1,1,1,1,1)+u(-1,-1,-1,-1,-1,-1,-1,-1)\\
&&-u(1,1,-1,-1,1,1,-1,-1)-u(-1,-1,1,1,-1,-1,1,1)\\
&&+u(1,1,-1,-1,-1,-1,1,1,)+u(-1,-1,1,1,1,1,-1,-1)\\
&& -u(1,-1,1,-1,1,-1,1,-1)-u(-1,1,-1,1,1,-1,1,-1)\\
&&-u(1,-1,1,-1,-1,1,-1,1)-u(-1,1,-1,1,-1,1,-1,1)\\
&&-u(1,-1,-1,1,1,-1,-1,1) -u(-1,1,1,-1,-1,1,1,-1)\\
&& +u(1,-1,-1,1,-1,1,1,-1)+u(-1,1,1,-1,1,-1,-1,1)\\
& \in & \gD_{8,8}^+.
\eean
According to (\ref{sk2}) we obtain $\,\langle \eta,\eta \rangle =16$.

\subsection{$H=\Sp(p,\CC)\cdot \SL(2,\CC)\;\;\; (r=4p,\, s=4p)$}

The Lie algebra $\,\fsp(p) \oplus \fsp(1)\,$ is the compact real form of 
$\,\fsp(p,\CC) \oplus \fsl(2,\CC)\,$. Hence, $(\fsp(p) \oplus \fsp(1)) \oplus (\fsp(p) \oplus 
\fsp(1)) \subset \fso(4p) \oplus \fso(4p) \,$ is the compact real form of 
$\,(i_{\RR^*}(\fsp(p,\CC) \oplus \fsl(2,\CC)))^{{\mathbb C}}\,$ and we have 
\[
\mbox{dim}\,V_{\fsp(p,\CC)\oplus \fsl(2,\CC)} = 
\,\,\mbox{dim}\,V_{(\fsp(p) \oplus \fsp(1)) \oplus (\fsp(p) \oplus \fsp(1))}. 
\]
Because of 
\[
\gD_{8p}|_{(\fsp(p) \oplus \fsp(1)) \oplus (\fsp(p) \oplus \fsp(1))}\, \cong 
\, \gD_{4p}|_{\fsp(p) \oplus \fsp(1)} \otimes \, \gD_{4p}|_{\fsp(p) \oplus 
\fsp(1)} 
\]
and $\,V_{\fsp(p) \oplus \fsp(1)} = \{0\}\,$ we obtain $\,V_{\fsp(p,\CC) \oplus 
\fsl(2,\CC)} = \{0\}.$ 
\vspace{1cm}\\
Summarizing the previous calculations we obtain the theorem formulated in the 
introduction. 

{\small

\vspace{10mm}

\noindent
Helga Baum\\
Ines Kath\\
Institut f\"ur Mathematik\\
Humboldt-Universit\"at Berlin\\
Sitz: Ziegelstr. 13 a\\
10099 Berlin\\[1ex]
baum@mathematik.hu-berlin.de\\
kath@mathematik.hu-berlin.de
} 
\edo